\documentclass[11pt]{article}

\usepackage{amscd,amsmath, amssymb, fancyhdr, url}

\usepackage{fourier, libertine}

\newcommand{\version}{version 1.0,\ \ January 27, 2016}
\newcommand{\f}{\varphi}

\newcommand{\RR}{\mathbb{R}}

\numberwithin{equation}{section}

\def\eqref#1{(\ref{#1})}
\newcommand{\goth}{\mathfrak}

\newcommand{\Z}{{\mathbb Z}}
\newcommand{\C}{{\mathbb C}}
\newcommand{\R}{{\mathbb R}}
\newcommand{\Q}{{\mathbb Q}}

\newcommand{\6}{\partial}
\def\1{\sqrt{-1}\:}

\newcommand{\cntrct}                
{\hspace{2pt}\raisebox{1pt}{\text{$\lrcorner$}}\hspace{2pt}}
\newcommand{\arrow}{{\:\longrightarrow\:}}

\newcommand{\cac}{{\cal C}}

\renewcommand{\bar}{\overline}
\renewcommand{\phi}{\varphi}
\renewcommand{\epsilon}{\varepsilon}
\renewcommand{\geq}{\geqslant}

\newcommand{\Aut}{\operatorname{Aut}}

\renewcommand{\Re}{\operatorname{Re}}

\newcounter{Mycounter}[section]
\newcounter{lemma}[section]
\newcounter{claim}[section]
\renewcommand{\theclaim}{{Claim \thesection.\arabic{claim}}}
\newcommand{\claim}{%
     \setcounter{claim}{\value{Mycounter}}
     \refstepcounter{claim}
     \stepcounter{Mycounter}
     {\noindent \bf \theclaim:\ }}

\newcounter{sublemma}[section]
\newcounter{corollary}[section]
\renewcommand{\thecorollary}{{Corollary \thesection.\arabic{corollary}}}
\newcommand{\corollary}{%
     \setcounter{corollary}{\value{Mycounter}}
     \refstepcounter{corollary}
     \stepcounter{Mycounter}
     {\noindent \bf \thecorollary:\ }}

\newcounter{theorem}[section]
\renewcommand{\thetheorem}{{Theorem \thesection.\arabic{theorem}}}
\newcommand{\theorem}{%
     \setcounter{theorem}{\value{Mycounter}}
     \refstepcounter{theorem}
     \stepcounter{Mycounter}
     {\noindent \bf \thetheorem:\ }}

\newcounter{conjecture}[section]
\newcounter{proposition}[section]
\renewcommand{\theproposition} {{Proposition \thesection.\arabic{proposition}}}
\newcommand{\proposition}{%
     \setcounter{proposition}{\value{Mycounter}}
     \refstepcounter{proposition}
     \stepcounter{Mycounter}
     {\noindent \bf \theproposition:\ }}

\newcounter{definition}[section]
\renewcommand{\thedefinition} {{Definition~\thesection.\arabic{definition}}}
\newcommand{\definition}{%
     \setcounter{definition}{\value{Mycounter}}
     \refstepcounter{definition}
     \stepcounter{Mycounter}
     {\noindent \bf \thedefinition:\ }}

\newcounter{example}[section]
\newcounter{remark}[section]
\renewcommand{\theremark}{{Remark \thesection.\arabic{remark}}}
\newcommand{\remark}{%
     \setcounter{remark}{\value{Mycounter}}
     \refstepcounter{remark}
     \stepcounter{Mycounter}
     {\noindent \bf \theremark:\ }}

\newcounter{problem}[section]
\newcounter{question}[section]
\makeatletter

\@addtoreset{equation}{section}
\@addtoreset{footnote}{section}
\makeatother

\def\blacksquare{\hbox{\vrule width 5pt height 5pt depth 0pt}}
\def\endproof{\blacksquare}

\begin{document}

\begin{center}
{\LARGE\bf
LCK rank of locally conformally K\"ahler manifolds with potential
}\\[3mm]
{\large
Liviu Ornea\footnote{Partially supported by CNCS UEFISCDI, project number PN-II-ID-PCE-2011-3-0118.}, and 
Misha
Verbitsky\footnote{Partially supported 
by RSCF grant 14-21-00053 within AG Laboratory NRU-HSE.\\[1mm]

}\\[4mm]
}

\end{center}


{\small
\hspace{0.15\linewidth}
\begin{minipage}[t]{0.7\linewidth}
{\bf Abstract} \\ 
An LCK manifold with potential is a compact quotient of a
K\"ahler manifold $X$ equipped with a positive K\"ahler potential
$f$, such that the monodromy group acts on $X$ by holomorphic homotheties
and multiplies $f$ by a character. The LCK rank
is the rank of the image of this character, considered
as a function from the monodromy group to real numbers. We
prove that an LCK manifold with potential can have any
rank between 1 and $b_1(M)$. Moreover, LCK manifolds
with proper potential (ones with rank 1) are dense.
Two {\em errata} to our previous work are given in the last Section.
\end{minipage}
}

\tableofcontents

\hfill

\hfill
\section{Introduction}

\subsection{LCK manifolds}

A complex manifold $(M,I)$ is called {\bf locally conformally K\"ahler} (LCK) if it
admits a  Hermitian metric $g$ and a closed 1-form $\theta$,  called {\bf the Lee
form}, such that the  
fundamental 2-form $\omega(\cdot,\cdot):=g(\cdot, I\cdot)$
satisfies the integrability condition
\begin{equation}\label{deflck}
d\omega=\theta\wedge\omega,\quad d\theta=0.
\end{equation}

The above definition  is equivalent (see \cite{DO}) to the existence
of a  covering $\tilde M$ endowed with a K\"ahler metric $\Omega$ which is
acted on by the deck group $\Aut_M(\tilde M)$ by holomorphic homotheties. Hence, if $\tau\in \Aut_M(\tilde M)$, then $\tau^*\Omega=c_\tau\cdot\Omega$, where $c_\tau\in \R^{>0}$ is the scale factor. This defines a character 
\begin{equation}\label{chi}
\chi:\Aut_M(\tilde M)\longrightarrow \R^{>0},\quad \chi(\tau)=c_\tau.
\end{equation}

Two subclasses of LCK manifolds will be of interest to us.

The {\bf Vaisman} class is formed by LCK manifold $(M,\omega, \theta)$ with parallel Lee form with respect  to the
Levi-Civita connection of $g$. While the LCK
condition is conformally invariant (if $g$ is LCK,
then any $e^f\cdot g$ is still LCK, with Lee form $\theta+df$), the Vaisman condition
is not. The main example of Vaisman manifold is the
diagonal Hopf manifold (\cite{_OV:Shells_}). Also, all compact complex submanifolds of Vasiman manifolds are Vaisman, too, \cite{_Verbitsky:LCHK_}. The Vaisman compact complex
surfaces are classified in \cite{be}.

We observed in \cite{_Verbitsky:LCHK_}, \cite{ov01} that the K\"ahler form of the universal cover of any Vaisman manifold has global potential represented by the square of the length 
of the Lee form. Moreover, the deck group acts on the potential by multiplying it with the character $\chi$. This led us to introducing the larger class of LCK manifolds 
 {\bf with potential}. The precise definition requires the 
 existence 
of a K\"ahler covering on which the K\"ahler metric
has global, positive and proper potential function 
which
is acted on by homotheties by the deck group. Besides Vaisman manifolds, there exist non-Vaisman
examples, such as the non-diagonal Hopf manifolds, \cite{ov01}. 

\subsection{LCK manifolds with potential}

``LCK manifolds with potential'' can be defined as 
LCK manifolds $(M, \omega, \theta)$ equipped with a smooth function
$\psi\in \mathcal{C}^\infty (M)$,
\begin{equation}\label{_pote_impro_Equation_}
\omega= d_\theta d^c_\theta \psi,
\end{equation}
where $d_\theta(x)=dx-\theta\wedge x$, $d_\theta^c= I d_\theta I^{-1}$,
and the following properties are satisfied:
\begin{equation}\label{_pote_pro_restrictions_Equation_}
\begin{minipage}[m]{0.92\linewidth}
\begin{description}
\item[(i)] $\psi >0$;
\item[(ii)] the  class $[\theta]\in H^1(M, \R)$ is proportional to a rational one.
\end{description}
\end{minipage}
\end{equation}
For more details and historical context of this definition,
please see Subsection \ref{_LCK-Pot-Subsection_}.
The differential $d_\theta$ is identified
with the de Rham differential with coefficients in a flat
line bundle $L$ called {\bf the weight bundle}.
In this context, $\psi$ should be considered as a section of $L$.
After passing to the smallest covering $\tilde M\stackrel \pi \arrow M$ where $\theta$
becomes exact, the pull-back bundle $\pi^* L$ can be trivialized
by a parallel section. Then the equation \eqref{_pote_impro_Equation_}
becomes $\tilde \omega= dd^c\tilde  \psi$, where $\tilde \omega$
is a K\"ahler form on $\tilde M$, and $\tilde \psi$ the 
K\"ahler potential.

Since $\tilde M\stackrel \pi \arrow M$ is the smallest covering
where $\theta$ becomes exact, its monodromy
is equal to $\Z^k$, where $k$ is the rank of the smallest rational
subspace $V\subset H^1(M, \Q)$ such that $V\otimes_\Q \R$
contain $[\theta]$.
In particular, the condition \eqref{_pote_pro_restrictions_Equation_} (ii)
means precisely that $\tilde M\stackrel \pi \arrow M$ is a $\Z$-covering.
This implies that the definition 
\eqref{_pote_impro_Equation_}-\eqref{_pote_pro_restrictions_Equation_} 
is equivalent to the historical one (\ref{_LCK_w_p_orig_Definition_}).

However, the condition
\eqref{_pote_pro_restrictions_Equation_} (ii) 
is more complicated: there are examples of
LCK manifolds satisfying \eqref{_pote_impro_Equation_}
and not \eqref{_pote_pro_restrictions_Equation_} (ii)
(Subsection \ref{_impro_pote_Subsection_}).
Still, any complex manifold $(M,\omega, \theta)$ admitting
an LCK metric   with potential $\psi$
satisfying  \eqref{_pote_impro_Equation_}, admits an
LCK metric satisfying 
\eqref{_pote_impro_Equation_}-\eqref{_pote_pro_restrictions_Equation_}
in any $\mathcal{C}^\infty$-neighbourhood of $(\omega, \theta)$.
Therefore the condition
\eqref{_pote_pro_restrictions_Equation_} (ii)
is not restrictive, and for most applications,
unnecessary.

It makes sense to modify the notion of LCK manifold 
with potential to include the following notion
(Subsection \ref{_impro_pote_Subsection_}):

\hfill

\definition\label{_lck_proper_improper_Definition_}
Let $(M, \omega, \theta)$ be an LCK manifold, and $\psi\in
\mathcal{C}^\infty (M)$ a positive function satisfying $d_\theta d^c_\theta \psi=\omega$. 
Denote by $k$ the rank of the smallest rational
subspace $V\subset H^1(M, \Q)$ such that $V\otimes_\Q \R$
contain $[\theta]$. Then $\psi$ is called {\bf proper potential}
if $k=1$ and {\bf improper potential} if $k >1$.

\subsection{Some errors found}

This paper is much influenced by Paul Gauduchon,
who discovered an error in our result mentioned as obvious in
\cite{ov_imrn_10}. In \cite{ov_imrn_10}, we claimed erroneously
that an LCK metric is pluricanonical, see \cite{ko}, if and only if 
it admits an LCK potential. This was obvious because
(as we claimed) the equations for LCK with potential
and for pluricanonical metric are the same.
Unfortunately, a scalar multiplier was missing
in our equation for the pluricanonical metric.

From an attempt to understand what is brought by the
missing multiplier, this paper grew, and we found
an even stronger result: any compact pluricanonical manifold
is Vaisman. Very recently, Andrei and Sergiu Moroianu gave a simple, direct proof of this result, using elegant tensor computations, \cite{mm}.

However, during our work trying to plug a seemingly harmless
mistake, we discovered a much more offensive error, which 
has proliferated in a number of our papers.

In \cite{_OV:_Structure_}, we claimed that any Vaisman manifold
admits a $\Z$-covering which is K\"ahler. This is true for
locally conformally hyperk\"ahler manifolds, as shown in 
\cite{_Verbitsky:LCHK_}. However, this result is false
for more general Vaisman manifolds, such as a Kodaira surface
(\ref{_Lee_open_Theorem_}).

It is easiest to state this problem and its solution
using the notion of ``LCK rank'' (\ref{_LCK_rank_Definition_}),
defined in \cite{gopp}  and studied in \cite{PV}.
Briefly, the LCK rank is the smallest $r$ such that 
there exists
a $\Z^r$-covering $\tilde M$ of $M$ such that the pullback
of the LCK metric is conformally equivalent
to a K\"ahler metric on $\tilde M$.

It turns out that the LCK rank of
a Vaisman manifold can be any number between
1 and $b_1(M)$ (\ref{_Lee_open_Theorem_}).
Moreover, for each $r$, the set of all Vaisman
metrics of LCK rank $r$ is dense in the
space of all Vaisman metrics (say, with $\mathcal{C}^\infty$-topology).

It is disappointing to us (and even somewhat alarming) that 
nobody has discovered this important error earlier.

However, not much is lost, because 
the metrics which satisfy the Structure Theorem of 
\cite{_OV:_Structure_} are dense in the space of all LCK metrics, hence all results
of complex analytic nature remain true. To 
make the remaining ones correct, we need to add
``Vaisman manifold of LCK rank 1'' or
``Vaisman manifold with proper potential'' (Subsection 
\ref{_impro_pote_Subsection_}) to the set of assumptions
whenever \cite{_OV:_Structure_} is used. 

Still, we want to offer our apologies to the
mathematical community for managing to mislead our colleagues
for such a long time.

For more details about our error
and an explanation where the arguments of 
\cite{_OV:_Structure_}  failed, 
please see Subsection \ref{_LCK_rank_errata_Subsection_}.


\section{LCK manifolds: properness  of the potential}
\label{_LCK-Pot-proper_Section_}

\subsection{LCK manifolds with potential: historical
  definition}
\label{_LCK-Pot-Subsection_}

When the notion of LCK manifold with potential was introduced in 
\cite{ov01}, we assumed properness of the potential. Later, it was ``proven''
that the potential is always proper
(\cite{_OV_Automorphisms_}). Unfortunately, the proof was false (see the  {\em  Errata} to this paper,
Section \ref{err}). In  view of this error and other results in
Section \ref{err}, it makes sense to generalize the notion of
LCK manifold with potential to include the manifolds with 
LCK rank $>1$. For the old notion of LCK with potential
we should attach ``proper'' to signify that the potential
is a proper function on the minimal K\"ahler covering.

\hfill

\definition \label{_LCK_w_p_orig_Definition_}
(\cite{ov01}) 
An {\bf LCK manifold with proper potential} is a manifold
which admits a K\"ahler covering $(\tilde M, \tilde \omega)$ and a 
smooth function $\phi:\,\tilde M \rightarrow \RR^{>0}$ 
(the {\bf LCK potential}) satisfying the following conditions:
\begin{description}
\item[(i)]  $\phi$ is proper, \emph{i.e.} its level sets are compact;
\item[(ii)] The deck transform group acts on $\phi$ by multiplication 
with the  character $\chi$ (see \eqref{chi}): 
$\tau^* (\phi)=\chi(\tau) \phi$, where $\tau\in \Aut_M(\tilde M)$ 
is any deck transform map\footnote{In general, differential forms $\eta\in \Lambda^\bullet \tilde M$  which satisfy $\tau^* \eta=\chi(\tau) \eta$ are  called {\bf automorphic}. In particular, so is the K\"ahler form on $\tilde M$.}.
\item[(iii)] $\phi$ is a  K\"ahler potential, \emph{i.e.} 
$dd^c\phi = \tilde \omega$.
\end{description}
 
\hfill

\remark
In this situation, we can fix a choice of 
LCK metric on $M$ by writing $\pi^*\omega= \phi^{-1}dd^c\phi$.
Further on, we shall tacitly assume that this choice is used
whenever we work with an LCK manifold with potential.
In this case, the Lee form is written as $\pi^*\theta=d\log\phi$,
and $d^c(\pi^*\theta)=\pi^*\theta\wedge I(\pi^*\theta)-\pi^*\omega$ 
(\cite{_OV:MN_}).

\hfill

\remark Positivity of the potential cannot be relaxed, as the following simple example (for which we thank V. Vuletescu) shows. On $\C^2\setminus 0$, with $\Z$ acting as $(z_1, z_2)\mapsto (2z_1, 2z_2)$ (the quotient being the usual Hopf surface, which is Vaisman) take
$\phi(z_1,z_2)=|z_1|^2 + |z_2|^2 - \frac{1}{3}(z_1+\bar z_1)^2$. 
Then:
$\6\bar \6 \phi =\frac{1}{3}dz_1\wedge d\bar z_1 + dz_2\wedge d\bar z_2$ 
and
$\phi(2z_1, 2z_2)=4\phi(z_1, z_2)$, 
and hence the potential is automorphic.
But 
$\phi(1,0)=-\frac{1}{3}$, 
$\phi(0,1)=1$, and $\phi^{-1}(0)$ is non-empty.

More general examples are obtained by starting from any automorphic potential
and adding the real part of a convenient holomorphic function, automorphic with the same automorphy factor as the potential.

\subsection{Properness of the LCK potential}

In \cite{ov01}, it was also shown that the properness condition 
is equivalent to the following condition on
the deck transform group of $\tilde M$. Recall
that a group is {\bf virtually cyclic} if it contains
$\Z$ as a finite index subgroup. We obtain the following claim (which proof we include for convenience):

\hfill

\claim\label{_proper_virt_cy_Claim_}
Let $M$ be a compact manifold, 
$\tilde M$ a covering, and 
$\phi:\,\tilde M \rightarrow \R^{>0}$ 
an automorphic function, that is, a
function which satisfies $\gamma^* \phi= c_\gamma \phi$
for any deck transform map $\gamma$, where $c_\gamma=\chi(\gamma)\neq 1$ 
is constant ($\chi$ being the character \eqref{chi}).
Then $\phi$ is proper if and only if the deck transform group $\Gamma:=\Aut_M(\tilde M)$ of
$\tilde M$ is virtually cyclic. 

\hfill

{\bf Proof:} By passing to a smaller cover, we may suppose $\Gamma=\Z$ and 
let $\gamma$ be a generator of $\Z$, such that
$\gamma^*\phi=\lambda\phi$, and $\pi:\; \tilde M \arrow
\tilde M/\Gamma =M$ the quotient map.
{ Then $\phi^{-1}([1, \lambda[)$ is a fundamental domain
	of the $\Gamma$-action.} Therefore
$\pi:\;\phi^{-1}([1,\sqrt  \lambda])\arrow M$ is bijective
onto its image, which is compact, and hence 
$\phi^{-1}([1,\sqrt  \lambda])$ is also compact.
This implies that 
{the  preimage of any closed interval is compact}.

Conversely, suppose $\phi$ is proper and, by absurd,  
assume $\Gamma\neq \Z$. Then { 
	$\Gamma$ is a dense subgroup of $\R^{>0}$.}
Fix $x\in \tilde M$ and a nonempty interval
$]a,b[\subset \R^{>0}$, and let 
$${\goth H}:=\left\{\gamma\in \Gamma\,;\, \phi(\gamma(x))\in ]a, b[\right\}.$$ 
Since $\Gamma$ is dense,  ${\goth H}$ is infinite.
However, $\phi({\goth H}\cdot x)\subset [a, b]$, hence
{ the infinite discrete set ${\goth H}\cdot x$
	is contained in a compact $\phi^{-1}([a,b])$.}
This contradiction ends the proof. \endproof

\hfill

\definition \label{_LCK_rank_Definition_}
Let $(M,\omega,\theta)$ be an LCK manifold. Define the 
{\bf LCK rank} as the dimension of the smallest
rational subspace $V\subset H^1(M,\Q)$ such that
the Lee class $[\theta]$ lies in $V\otimes_\Q \R$.

\hfill

\remark
The character $\chi:\; \Aut_M(\tilde M)\arrow \R^{>0}$
is defined on any LCK manifold, because the  K\"ahler form $\tilde \omega$ is automorphic by definition:
$\tau^* (\tilde \omega)=\chi(\tau) \tilde \omega$. Then one can see that the LCK rank as defined above coincides with the rank of 
the image of $\chi:\; \Aut_M(\tilde M)\arrow \R^{>0}$ which is also called 
{\bf the weight monodromy group} of the  LCK manifold. See also \cite{gopp} for another interpretation of the LCK rank  and see \cite{PV} for examples on non-Vaisman compact LCK manifolds with K\"ahler rank greater than 1. Clearly, LCK rank 0 corresponds to globally conformally K\"ahler structures.

From \ref{_proper_virt_cy_Claim_} above 
it follows that condition (i) in \ref{_LCK_w_p_orig_Definition_} is equivalent to
$M$ being of LCK rank 1. 

\hfill

In \cite{_OV:MN_}, we managed to get rid of the 
need to take the covering in \ref{_LCK_w_p_orig_Definition_},
by using the Morse-Novikov (twisted) differential $d_\theta:=d-\theta\wedge\cdot$,
where $\theta\wedge\cdot(x)=\theta\wedge x$, and $\theta$
is the Lee form. In \cite{_OV:MN_} the definition
of LCK manifold with potential 
was restated equivalently as follows.


\hfill

\definition\label{_LCK_pot_MN_Definition_}
Let $(M,\omega,\theta)$ be an LCK manifold of LCK rank 1.
Then $M$ is called {\bf LCK manifold with potential}
if there exists a positive function $\phi_0\in \mathcal{C}^\infty (M)$ satisfying 
$d_\theta d^c_\theta(\phi_0)=\omega$, where $d^c_\theta =
I d_\theta I^{-1}$.

\hfill

\claim\label{_MN_def_equi_Claim_}
 \ref{_LCK_w_p_orig_Definition_} is equivalent to
\ref{_LCK_pot_MN_Definition_}.

\hfill

{\bf Proof:}
To see that \ref{_LCK_w_p_orig_Definition_} and
\ref{_LCK_pot_MN_Definition_} are equivalent, consider
the smallest covering $\pi:\; \tilde M \arrow M$ such that
$\pi^*\theta$ is exact, and take a function $\psi$
satisfying $d\psi=\pi^*\theta$. Since $\pi^*\theta$ is invariant
under the deck transform group $\Gamma$,
for each $\gamma\in \Gamma$ one has $\gamma^*\psi=\psi +
c_\gamma$, where $c_\gamma $ is a constant. Consider the
multiplicative character $\chi:\; \Gamma\arrow \R^{>0}$
given by $\chi(\gamma)= e^{c_\gamma}$. Let
$\Lambda_\chi^*(M)$ denote the space 
of automorphic forms on $\tilde M$
which satisfy $\gamma^*\eta= \chi_\gamma\Gamma\eta$.
The map $\Lambda^*(M) \stackrel \Psi \arrow \Lambda_\chi^*(M)$
mapping $\eta$ to $e^{-\psi}\pi^*\eta$ makes the following
diagram commutative:
\[
\begin{CD}
\Lambda^*(M) @>\Psi>> \Lambda_\chi^*(M)\\
@V {d_\theta} VV @V {d} VV\\
 \Lambda^*(M) @>\Psi>> \Lambda_\chi^*(M)
\end{CD}
\]
Then $\Psi$ maps a ``potential'' $\phi_0$ in the sense of 
\ref{_LCK_pot_MN_Definition_} to a potential $\psi$
in the sense of \ref{_LCK_w_p_orig_Definition_} and vice versa.
Properness of $\Psi(\phi_0)$ is equivalent to 
$\Gamma$ being virtually cyclic, as 
\ref{_proper_virt_cy_Claim_} implies.
The existence of a K\"ahler covering with 
virtually cyclic deck transform group is
clearly equivalent to $M$ having LCK rank 1. \endproof

\hfill 

\begin{example}
On a diagonal Hopf manifold $(\C^n\setminus 0)/\Z$, with
LCK form and Lee form written on $\C^n\setminus 0$
respectively  $\omega=|z|^{-2}\sum dz^i\wedge d\bar z^i$
and $\theta=-d\log|z^2|$, the function $\f_0\in\cac(M)$ is
the constant function 1,
while the potential function on $\C^n\setminus 0$ is $|z^2|$.
\end{example}

\hfill

We now show that automorphic potentials can be
approximated by proper ones. The following argument is
taken from \cite{_OV:MN_}.

\hfill

\claim\label{_LCK_appro_Claim_}
Let $(M,\omega, \theta)$ be an LCK manifold, and
$\phi\in \mathcal{C}^\infty (M)$ a function satisfying 
$d_\theta d^c_\theta(\phi)=\omega$.
Then $M$ admits an LCK structure $(\omega', \theta')$
of LCK rank 1, approximating $(\omega, \theta)$ in $\mathcal{C}^\infty$-topology.

\hfill

{\bf Proof:} Replace $\theta$ by a form $\theta'$ with rational
cohomology class $[\theta']$ in a sufficiently small
$\mathcal{C}^\infty$-neighbourhood of $\theta$, and 
let $\omega':=d_{\theta'} d_{\theta'}^c(\phi)$.
Then $\omega'$ approximates $\omega$
in $\mathcal{C}^\infty$-topology, hence for $\theta'$
sufficiently close to $\theta$, the form $\omega'$
is positive. It is $d_{\theta'}$-closed, because
$d_{\theta'}^2=0$, hence $0=d_{\theta'}\omega'=d\omega'-\theta'\wedge\omega'$.  
This implies that $(\omega', \theta')$
is an LCK structure. The K\"ahler rank of
an LCK manifold is the dimension of the smallest
rational subspace $W\subset H^1(M, \Q)$ such that
$W\otimes_\Q \R$ contains the cohomology class of the
Lee form. Since $[\theta']$
is rational, $(M,\omega', \theta')$ has LCK rank 1. \endproof

\subsection{LCK manifolds with proper and improper
  potential}
\label{_impro_pote_Subsection_}

It seems now that the equation 
$d_\theta d^c_\theta \psi=\omega$ (on the LCK manifold $M$ itself) 
is more fundamental than the notion of LCK manifold with 
(proper) potential. For most applications, this (more general)
condition is already sufficient. 

The relation between manifolds with 
$d_\theta d^c_\theta \psi=\omega$
and LCK with potential is similar to the relation
between general Vaisman manifolds and quasiregular ones\footnote{(Quasi)Regularity and irregularity of a Vaisman manifold refers to the (quasi)regularity and irregularity of the 2-dimensional canonical foliation generated by $\theta^\sharp$ and $I\theta^\sharp$.}.
One could always deform an irregular Vaisman manifold
to a quasiregular one, and quasiregular Vaisman manifolds
are dense in the space of all Vaisman manifolds.

The notion of ``LCK manifold with improper potential''
is similar, in this regard, to the notion of
irregular Vaisman or irregular Sasakian manifold\footnote{Here (ir)regularity refers to the 1-dimensional foliation generated by the Reeb field.}, \cite{bg}.

\hfill

\definition
Let $(M, \theta, \omega)$ be an LCK manifold,
and $\psi$ a strictly positive function 
which satisfies $d_\theta d^c_\theta\psi=\omega$.
Then $(M, \theta, \omega)$ is called {\bf a manifold
with improper LCK potential} if its LCK rank is $\geq 2$,
and {\bf a manifold with proper LCK potential}
if it has LCK rank 1.

\hfill

\remark The expressions ``proper potential'' and
``improper potential'', when used  for solutions of the
equation $d_\theta d^c_\theta\psi=\omega$, as in the above
Definition, do not refer to the properness of 
$\psi:M\longrightarrow \R$, which is not plurisubharmonic and is always proper if $M$ is
compact.

\hfill

Note that ``LCK with potential''
was previously used instead of ``manifold with proper LCK
potential''; now (in light of the discovery of Vaisman
manifolds having improper potential, see Section
\ref{err})  it makes sense to change the terminology
by including improper potentials in the definition
of LCK with potential.

\hfill

\ref{_LCK_appro_Claim_} can be rephrased as follows.

\hfill

\proposition
Let $(M, \omega, \theta, \psi)$ be a compact LCK manifold
with improper LCK potential. Then $(\omega, \theta, \psi)$
can be approximated in the $\mathcal{C}^\infty$-topology by an LCK
structure with proper LCK potential.
\endproof

\hfill

\remark We have just proven that existence of an LCK metric with 
improper LCK potential implies existence of a metric with
proper LCK potential on the same manifold. The converse is
clearly false: when $H^1(M, \R)$ is 1-dimensional, any
Lee class is proportional to an integral cohomology class,
and any LCK structure has LCK rank 1, hence $M$ admits 
no metrics with improper LCK potentials. 

However, in all other situations improper potentials do
exist.

\hfill

\proposition\label{_impro_pote_exists_Proposition_}
Let $(M, \omega, \theta, \psi)$ be an LCK manifold with potential, and suppose 
$b_1(M)>1$. Then $M$ admits an
LCK metric $(M, \omega', \theta', \psi)$ with improper
potential and arbitrary LCK rank between 2 and $b_1(M)$.
Moreover, $(\omega', \theta')$ can be chosen in arbitrary
$\mathcal{C}^\infty$-neighbourhood of $(\omega, \theta)$.

\hfill

{\bf Proof:} Choose a closed $\theta'$ in a sufficiently
small neighbourhood of $\theta$, and let $V_\theta$ be the 
smallest rational subspace of $H^1(M, \R)$ such that
$V_\theta\otimes_\Q \R$ contains $\theta$. Since the
choice of the cohomology class $[\theta']$ is arbitrary
in a neighbourhood of $[\theta]$, the dimension of
$V_\theta$ can be chosen in arbitrary way. Choosing 
$\theta'$ sufficiently close to $\theta$, we can
assume that the (1,1)-form $\omega':=d_{\theta'}
d^c_{\theta'}(\psi)$ is positive definite. 
Then $(M, \omega', \theta', \psi)$ is an LCK manifold
with improper potential and arbitrary LCK rank.
\endproof

\section{\em Errata}\label{err}


\subsection{Pluricanonical condition revisited}

 In Section 3 of \cite{ov_imrn_10} the following erroneous
claim was made: ``We now prove that the pluricanonical condition is
equivalent with the existence of an automorphic potential
on a K\"ahler covering.''

Then we proceeded to make calculations purporting to show
that pluricanonical condition is equivalent to the LCK with
potential condition $d(I\theta) = \omega - \theta\wedge I\theta.$
Here, the scalar term is lost: the correct
equation (in the notation  of \cite{ov_imrn_10})
is $d(I\theta) = |\theta|^2\omega - \theta\wedge I\theta.$

It is of course true that this equation, indeed,
implies $d(I\theta) = \omega - \theta\wedge I\theta$.

However, the converse statement is false:
not all LCK manifolds
with potential admit a pluricanonical LCK structure, but
only Vaisman ones, see \cite{mm}.

\subsection{LCK rank of Vaisman manifolds}
\label{_LCK_rank_errata_Subsection_}

Recall that the {\bf LCK rank} of an LCK manifold $(M, \omega,\theta)$ is
the rank of the smallest rational subspace $V$ in $H^1(M,\R)$
such that $V\otimes_\Q \R$ contains the cohomology class
$[\theta]$. When the LCK rank is 1, the manifold admits
a $\Z$-covering which is K\"ahler (Section \ref{_LCK-Pot-proper_Section_}).

In several papers published previously (\cite{_OV:_Structure_}, 
\cite{_OV:MN_}, \cite{_OV_Automorphisms_})
we claimed that a Vaisman manifold and an LCK manifold with potential
always have LCK rank 1. This is in fact false. In this section we 
produce a counterexample to these claims, and explain the error.

Notice, however, that, as we prove below,  any complex manifold which admits a
structure of a Vaisman manifold (or LCK manifold with potential)
{\em also} admits a structure of a Vaisman manifold (or LCK manifold with potential)
with LCK rank one. This means that all problems arising because of this error are of 
differential-geometrical nature; results of complex geometry remain valid.
This is probably the reason why the error was not noticed for so many years.
Moreover, {\em the set of Vaisman (or LCK with potential) structures with LCK rank 1
on a given manifold is dense in the set of all Vaisman (or LCK with potential)
structures.}

LCK manifolds with potential can have arbitrary
LCK rank, as follows from \ref{_impro_pote_exists_Proposition_}.
To construct a Vaisman manifold with an LCK rank 
bigger than 1, we proceed as follows. 

Let us recall some facts from Vaisman geometry used in this construction.
Any Vaisman manifold is equipped with a
canonical holomorphic foliation $\Sigma$, generated by
the Lee field $\theta^\sharp$ and $I(\theta^\sharp)$
(\cite{vai}, \cite{_Tsukada:foliation_}). This foliation
might have a global leaf space (in this case the Vaisman
manifold is called {\bf quasiregular}), or have non-closed
leaves ({\bf irregular Vaisman}). Locally, the leaf space always exists.
{\bf Transversal forms} are forms which are lifted (locally) from
the leaf space of $\Sigma$. K. Tsukada in \cite{_Tsukada:decomposition_} 
proved the following decomposition theorem.

\hfill

\theorem
The space of harmonic forms on a compact Vaisman manifold $M$
can be expressed as
\[
{\cal H}^*(M)= \theta \wedge {\cal H}^*_{tr}(M) \oplus {\cal H}^*_{tr}(M)
\]
where ${\cal H}^*_{tr}(M)$ is the space of transversal
harmonic forms.
\endproof

\hfill

The Vaisman manifold is transversally K\"ahler, that is,
the leaf space of the canonical foliation $\Sigma$ is
locally equipped with a complex structure and a 
globally defined transversally K\"ahler form.
This allows us to use the Hodge decomposition theorem
for transversal harmonic forms (\cite{elkacimi}),
entirely similar to the usual Hodge decomposition theorem
in K\"ahler geometry.
In particular, any transversal harmonic
1-form on a Vaisman manifold is the sum of a 
transversal holomorphic form and a transversal
antiholomorphic form.

This leads to the following useful corollary.

\hfill

\corollary\label{desc}
Let $M$ be a compact Vaisman manifold. Then the space
of harmonic 1-forms can be decomposed as
${\cal H}^1(M,\R)= \langle\theta\rangle \oplus \Re (H^0(\Omega^1_{tr}(M)))$,
where $\Omega^1_{tr}(M)$ denotes the sheaf of holomorphic
transversal 1-forms.
\endproof

\hfill

\proposition\label{_Vaisman_defo_transve_Proposition_}
Let $(M,\theta,\omega)$ be a compact Vaisman manifold,
$\alpha$ a harmonic 1-form, and $\theta':= \theta+\alpha$.
Consider the (1,1)-form $\omega':=d_{\theta'} d^c_{\theta'}(1)$ obtained as a
deformation of $\omega= d_\theta d^c_\theta(1)$. Assume
that $\alpha$ is chosen sufficiently small in such a way
that $\omega'$ is Hermitian (\ref{_impro_pote_exists_Proposition_}). 
Then $\omega'$ is conformally equivalent to a Vaisman form.

\hfill

{\bf Proof:}
Consider the holomorphic flow $F$ generated by the Lee field
$\theta^\sharp$. It fixes $\omega$ and $\theta$ and its lift $\tilde F$ to the universal cover  $\tilde M$ acts by non-trivial homotheties with respect to the K\"ahler metric $\tilde\omega$.

We shall show that: (1) $F$ preserves $\omega'$ and (2)  $\tilde F$ acts by non-trivial
homotheties with respect to the K\"ahler metric $\tilde\omega'$ corresponding to $\omega'$. 

As $\alpha$ is the sum of $\lambda\theta$ 
and a transversal form (\ref{desc}), it is
preserved by $F$, too. Then $\theta'$ is preserved by $F$,
and also $I(\theta')$ is preserved, because $F$ is
holomorphic. As $d$ commutes with the action of the flow
and $1$ is a constant function, $\omega'$ is preserved by
$F$. This proves (1).

As for (2), if $\tilde \omega$ is the K\"ahler form on $\tilde M$ corresponding to $(M,\omega)$, then note that $(\tilde M/\langle e^{\tilde F}\rangle$ is Vaisman, too, and hence has odd $b_1$, \cite{vai}. If $\tilde F$ acts by isometries on $\tilde \omega'$, then this K\"ahler metric descends to a K\"ahler metric on $\tilde M/\langle e^{\tilde F}\rangle$, contradiction.
 
By  \cite[Theorem A]{kor}, (1) and (2) 
imply that $\omega'$ is conformally equivalent
to a Vaisman metric.
\endproof

\hfill

This gives the following unexpected result:.

\hfill

\theorem\label{_Lee_open_Theorem_}
Let $M$ be a compact Vaisman manifold or LCK 
manifold with potential and let $\mathcal{L}\subset H^1(M, \R)$ be 
the set of cohomology classes of all Lee forms for the Vaisman
(LCK with potential) structures on $M$. Then $\mathcal{L}$ is open in $H^1(M,\R)$.

\hfill

{\bf Proof:} This is
\ref{_Vaisman_defo_transve_Proposition_} for
Vaisman manifolds and
\ref{_impro_pote_exists_Proposition_}
for LCK manifolds with potential.
\endproof

\hfill

We have shown that for a general Vaisman structure $(M, \omega, \theta)$, 
its LCK rank is equal to $b_1(M)$, and any number between 1
and $b_1(M)$ can be obtained as an LCK rank for an appropriate
choice of $\theta$.

\ref{_Lee_open_Theorem_} has the following consequences.

\hfill

\corollary
Let $M$ be a compact complex manifold which admits a
structure of a Vaisman manifold (or LCK manifold with potential)
$(M, \omega, \theta)$. Then $M$ admits a structure of 
a Vaisman manifold (or LCK manifold with potential) $(M, \omega', \theta')$ with
proper potential, that is, of LCK rank one. Moreover,
such $\omega'$ and $\theta'$ can be chosen in any 
neighbourhood of $(M, \omega, \theta)$.
\endproof

\hfill

Now, let us explain where the proof of \cite{_OV:_Structure_} 
(later refined in \cite{_OV_Automorphisms_}) failed. 

Let $M$ be a compact Vaisman manifold, and $\theta^\sharp$ its Lee
field. Then $\theta^\sharp$ 
acts on $M$ by holomorphic isometries, and on its smallest K\"ahler
covering $(\tilde M, \tilde \omega)$ by holomorphic homotheties. 
Denote by $G$ the closure of the group generated by $e^{t\theta^\sharp}$.
This group is a compact Lie group, because isometries form a compact
Lie group on a compact Riemannian manifold, and a closed subgroup of a Lie group is a Lie group
by Cartan's theorem. Moreover, it is commutative, because $\langle e^{t\theta^\sharp}\rangle$ is 
commutative, and this gives $G=(S^1)^{k}$.

Let $\tilde G$ be the group of pairs 
$(\tilde f\in \Aut(\tilde M), \ f\in G)$, making the following diagram
commutative:
\begin{equation*}
\begin{CD}
 \tilde M@>{\tilde f}>> \tilde M \\
@V{\pi}VV  @VV{\pi}V              \\
M@>{f}>>  M 
\end{CD}
\end{equation*}
Then $\tilde G$ is a covering of $G$, and the kernel of this projection
is $\tilde G\cap \Aut_{M}(\tilde M)$, where $\Aut_{\tilde M}(M)$ is the deck transform group of the
covering $\tilde M \arrow M$. 

Consider the homomorphism $\chi:\; \pi_1(M)\arrow \R^{>0}$
mapping an element of $\pi_1(M)$ considered as an automorphism
of $\tilde M$, to the K\"ahler homothety constant, 
$\gamma \mapsto \frac{\gamma^*\tilde \omega}{\tilde \omega}$.
Since $\tilde M$ is the smallest K\"ahler covering,
we identify $\Aut_{M}(\tilde M)$ with $\chi(\pi_1(M))\subset \R^{>0}$.

Now, let $\tilde G_0\subset \tilde G$ be the subgroup
acting on $\tilde M$ by isometries.
Since the group $\tilde G\cap \Aut_{M}(\tilde M)$ is a subgroup
of $\Aut_{\tilde M}(M)$, $\tilde G_0$ maps to its image in $G$ bijectively.

We assumed that $G_0$  (being the subgroup of elements of 
$\tilde G$ acting by isometries on both $\tilde M$ and $M$)
is closed in $G$. Then, if $\tilde G_0\cong S^{k-1}$,
this would imply that $\tilde G\cong (S^1)^{k-1} \times \R$,
proving that $M$ is a quotient of $\tilde M$ by a $\Z$-action.

However, this is false, because $G_0$ is closed in $\tilde G$,
but not closed in $G$. This is where the argument fails.


\hfill

\hfill

\noindent{\bf Acknowledgments.} L.O. thanks the Laboratory
for Algebraic Geometry at the Higher School of Economics
in Moscow for hospitality and excellent research
environment during February and April 2014, and April
2015.

Both authors are indebted to Paul 
Gauduchon and Andrei Moroianu, for drawing their attention
on the insufficient motivation of their arguments in
\cite{ov_imrn_10} and for illuminating discussions, and to
Victor Vuletescu for his advice, examples, and a very careful 
reading of the first draft of this paper.

{\small

\noindent {\sc Liviu Ornea\\
{\sc University of Bucharest, Faculty of Mathematics, \\14
Academiei str., 70109 Bucharest, Romania, \emph{and}\\
Institute of Mathematics ``Simion Stoilow" of the Romanian Academy,\\
21, Calea Grivitei Street
010702-Bucharest, Romania }\\
\tt liviu.ornea@imar.ro, \ \ lornea@fmi.unibuc.ro
}

\hfill

\noindent {\sc Misha Verbitsky\\
{\sc  Laboratory of Algebraic Geometry,
Faculty of Mathematics,\\ National Research University
Higher School of Economics,\\
7 Vavilova Str. Moscow, Russia, } \emph{and}\\
{\sc Universit\'e Libre de Bruxelles, D\'epartement de Math\'ematique\\
Campus de la Plaine, C.P. 218/01, Boulevard du Triomphe\\
B-1050 Brussels, Belgium\\
\tt verbit@mccme.ru
}}}

\end{document}